\documentclass[a4paper,12pt]{paper}
\usepackage{amsmath,latexsym,amssymb,newcom}
\usepackage[mathscr]{eucal}
\begin{document}
\title{\sf{On a time-discrete approach to solving
Navier-Stokes systems.}}
\author {\sf{Kazuhiro HORIHATA}}
\maketitle
\noindent{\bf DNS.} \rm \enspace
This article advocate a new approximate scheme for
 Navier-Stokes systems which is described below:
when we  let $\Omega$ and $S$ be a domain in $\mathbb{R}^N$ 
 and it's boundary $( N \, \in \, \mathbb{N} )$, $T$ a positive number
and put 
 $a \, \in \, \stackrel {\circ} {W_2^{1}} (\Omega ; \mathbb{R}^N)$
  with $\div a \, = \, 0$,
it is 
%
\begin{equation}
\left\{
\begin{array}{ll}
&
\dfrac {\partial  v}{\partial t}
 \, - \, \triangle v 
  \, + \, ( v \cdot D ) v 
   \, + \, \grad p \; = \; 0
\quad \mathrm{in} \; \Omega \times (0,T),
\\[3mm]
&
\div v \; = \; 0,
\\
&
\left. v \right\vert_{\Omega \times \{ 0 \}}
 \; = \; a
  \quad \mathrm{and} \; 
   \left. v \right\vert_{S \times (0,T)} \, = \, 0
\end{array}
\right.
\label{EQ:NS}
\end{equation}
for a function $v$ $:$ $\Omega \times (0,T)$
$\to$ $\mathbb{R}^N$.
\par
The symbols employed here follows from a book of \cite{ladyzhenskaya}.
\par
Our new scheme is to consider the following variational setting:
 Set $N_T$ $=$ $[T/h]$ and $n$ $\in$ $\{ 1,2, \cdots, N_T \}$, $h > 0$ 
  and put $v_0$ $=$ $a$:
Let us suppose a sequence of mappings $\{ v_k \}$ and 
 of functionals $\{ I_k \}$
  $( k \, = 1, \cdots , n-1)$ be given. 
   Then we adopt $v_n$ $( n = 1,2, \cdots, N_T)$ as a minimizer of the functional
\allowdisplaybreaks\begin{align}
& I_n [v] \; = \; \int_{\Omega}
 \biggl( \frac {|v(x) \, - v_{n-1} (x-hv_{n-1}(x))|^2}{2h} \, + \, 
  \frac 12 | D v (x)|^2 \biggr) \, dx
\label{DEF:VS}
\end{align}
in the class $\stackrel {\circ} {W_2^{1}} ( {\Omega} \, ; \, \mathbb{R}^{N})$
 $\oplus$ $\stackrel{\circ}{\mathbf{J}} (\Omega)$.
\par
Since $I [v]$ is convex and lower-semicontinuous, we can readily
 see that such a $v_n$ exists in $\stackrel {\circ} {W_2^{1}} ( {\Omega} \, ; \, \mathbb{R}^{N})$
  $\oplus$ $\stackrel{\circ}{\mathbf{J}} (\Omega)$.
\par
We hereafter call a sequence of $\{ v_n \}$ $( 1,2, \cdots, N_T)$ \textit{DNS}
 and note that $v_n$ belongs to 
  $W_2^{2}$ $( {\Omega} \, ; \, \mathbb{R}^{N})$
   and satisfies the Euler-Lagrange equation
\allowdisplaybreaks\begin{align}
&
\frac {v_n (x)\, - \, v_{n-1} (x-hv_{n-1}(x))} h
 \, - \, \triangle v_n (x)
  \, \in \, \stackrel{\circ}{\mathbf{J}^{\perp}} (\Omega).
\label{EQ:1st}
\end{align}
Reader should remark that the variational approach is not crucial to
 construct our \textit{DNS} 
  but important to discuss the nested term $v_{n-1} (x-hv_{n-1}(x))$
   inspired by the material derivative.
From Helmholtz-decomposition lemma, we state that there exists a function
 $p$ $\in$ $W_2^1 (\Omega)$ such that
\allowdisplaybreaks\begin{align}
&
\frac {v_n (x)\, - \, v_{n-1} (x-hv_{n-1}(x))} h
 \, - \, \triangle v_n (x)
  \; = \; D p (x),
\label{EQ:DNS}
\end{align}
that is equivalent to
\allowdisplaybreaks\begin{align}
&
\frac {v_n (x)\, - \, v_{n-1} (x)} h
\, + \, \int_0^1 
 v_{n-1} (x) \cdot D v_{n-1} (x - h \tau v_{n-1} (x)) \, d\tau
\notag
\\
&
\, - \, \triangle v_n (x)                                                                                                                                                        
 \; = \; D p (x)
\label{EQ:DNS-1}
\end{align}
in the sense of distribution.
This is a reason why we can regard \eqref{DEF:VS} or \eqref{EQ:1st} as an approximate formula
 for Navier-Stokes systems.
\par
Throughout the paper we assume that
\allowdisplaybreaks\begin{align}
&
| Dv_n | \; = \; O ( 1/\sqrt{h} ),
\tag*{(A)}
\label{EQ:Assume}
\end{align}
which seems quite natural because the differential coefficient of the
 first order with respect to {\it{time}} is as same as the twice
  {\it{space}}-differentials in the heat equation.

\noindent{\bf Results.} \rm \enspace
We introduce the energy decay estimate directly obtained from \eqref{EQ:1st}:
\begin{Thm}{\rm{(The Energy Estimate)}.} \label{THM:ED}
The following holds:
\allowdisplaybreaks\begin{align}
&
\frac h {2} \sum_{n=1}^{N_T}
 \int_{\Omega}
  \Bigl\vert 
   \frac {v_n (x) \, - \, v_{n-1} (x-hv_{n-1}(x))}{h}
    \Bigr\vert^2 \, dx
\, + \, \frac 12 
 \int_\Omega
  | D v_n (x) |^2 \, dx
\notag
\\
&
\; \le \;
 C e^{CT} \int_\Omega | D a (x) |^2 \, dx
\label{INEQ:DT-DVN}
\end{align}
for any positive integer $n$ in $\{ 1, 2, \cdots, N_T\}$,
 where $C$ is a positive constant independent of $h$.
\end{Thm}
{\underbar{{Proof of Theorem \ref{THM:ED}}}.}
\par
By substituting $v_{n-1} (x-hv_{n-1}(x)))$ as a comparative mapping for $v$ 
 in \eqref{DEF:VS}, extending it to $0$ outside $\Omega$,
  computing the change of variables $y$ $=$ $x \, - \, h v_{n-1} (x)$
   with $\div v_n$ $=$ $0$, and using assumption on $v_n$:
   \ref{EQ:Assume}, 
    we arrive at
\allowdisplaybreaks\begin{align}
&
\int_\Omega
 \frac {| v_n (x) \, - \, v_{n-1} (x)|^2} {2h}
  \, dx
\, + \,
 \int_\Omega | D v_n (x) |^2 \, dx
\label{INEQ:1}
\\
&
\; \le \; ( 1 + Ch ) 
 \int_\Omega | D v_{n-1} (x) |^2 \, dx,
\notag
\end{align}
where $C$ is a positive constant independent of $h$.
The recursive usage above enjoys \eqref{INEQ:DT-DVN}.
\qed
\par
In the below we prove the existence of 
 a weak solution of Navier-Stokes systems:
Before stating our theorem, we prepare a few symbols:
 Set $t_n$ $=$ $n h$ $(n=1,2,3,\cdots, N_T)$ and
$$ 
\aligned
v_{\bar h} (t,x) & \, \;= \, v_n (x) \\
 v_h (t,x) & \, \;= \, 
  \frac {t-t_{n-1}} h v_n (x)
   \, + \, 
    \frac {t_n - t} h v_{n-1} (x)
\endaligned
\aligned
\quad & \\
\quad & \quad (t_{n-1} < t \le t_n).
\endaligned
\label{DEF:Interpolate}
$$
\par
When no ambiguity may arise, we say a pair of functions
 $v_{\bar h}$ and $v_h$ to be \textit{DNS};
\begin{Thm}{\rm{(Main Theorem)}.} \label{THM:Main}
Under the hypothesis \ref{EQ:Assume},
{\it{DNS}} converge weakly-star to a function $v$ 
 in $V_2 (\Omega \times (0,T) ; \mathbb{R}^N)$
  as $h \searrow +0$.
Besides it converges strongly to $v$ in 
 $L^2 ( \Omega \times (0,T))$;
  $v$ is a weak solution of Navier-Stokes systems with
\allowdisplaybreaks\begin{align}
&
\int_0^T \, dt \int_\Omega
 \Bigl\langle v (x,t),
  \frac {\partial \phi}{\partial t} (x,t) \Bigr\rangle \, dx
\, - \, 
 \int_0^T \, dt \int_\Omega
  \langle D v (x,t), D \phi (x,t) \rangle \, dx
\; = \; 0
\label{EQ:NSWE}
\\
&
\mathrm{for} \; ^\forall \phi \, \in \, 
 \dot{C}^\infty ( \Omega \times (0,T)) 
  \; \mathrm{with} \; \div \phi \, = \, 0.
\notag
\end{align}
Furthermore $v$ satisifies
\begin{equation}
\int_\Omega | D v(x,t) |^2 \, dx
\; \le \; C e^{CT} \int_\Omega | D a (x)|^2 \, dx
\label{INEQ:NSED}
\end{equation}
for any time $t$ in $(0,T)$ where $C$ is a positive universal constant.
\end{Thm}
\underbar{{Proof of Theorem \ref{THM:Main}}}.
\par
The former is directly obtained by Theorem \ref{THM:ED} 
 combined with Rellich Kondrachev theorem and Poincar\'e inequality.
\par
Next we claim that the convergent function $v$ is actually a weak
 solution of Navier-Stokes systems;
Since \textit{DNS} implies
\allowdisplaybreaks\begin{align}
&
- \int_0^T \, dt \int_\Omega
 \Bigl\langle v_h (x,t),
  \frac {\partial \phi}{\partial t} (x,t) \Bigr\rangle \, dx
\, + \,
 \int_0^T \, dt \int_\Omega
  \langle D v_{\bar{h}} (x,t), D \phi (x,t) \rangle \, dx
\; = \; 0
\label{EQ:WDNS}
\\
&
\mathrm{for} \; ^\forall \phi \, \in \,
 \dot{C}^\infty (\Omega\times(0,T))
  \; \mathrm{with} \; \div \phi \, = \, 0,
\notag
\end{align}
using results in the former,
 we can pass to the limit of $h$ $\searrow$ $0$ to verify the latter.
%
\vskip 9pt
\bibliographystyle{amsalpha}

\end{document}